\documentclass{article}

\usepackage{amsmath}
\usepackage{amsfonts}
\usepackage{amssymb}

\newtheorem{theorem}{Theorem}

\newtheorem{lemma}[theorem]{Lemma}
\newtheorem{proposition}[theorem]{Proposition}

\newtheorem{problem}{Open Problem}
\newenvironment{proof}[1][Proof]{\textbf{#1.} }{\ $\square$}

\begin{document}

\title{Bryc's random fields: the existence and distributions analysis}
\author{
Wojciech Matysiak, Pawe{\l} J. Szab\l owski\\
Faculty of Mathematics and Information Science,\\
Warsaw University of Technology,\\
Pl. Politechniki 1,\\
00-661 Warsaw, Poland}

\date{}

\maketitle

\begin{abstract}
We examine problem of existence of stationary random fields with linear regressions and
quadratic conditional variances, introduced by Bryc \cite{bryc1}. Distributions of
the fields are identified and almost complete description of the possible sets of parameters
defining the first two conditional moments is given. This note almost solves
Bryc's problem concerning fields undetermined by moments - the only remaining set of parameters
for which the existence of Bryc's fields is unclear has Lebesgue measure zero.
\end{abstract}

\section{Introduction}

This note deals with square integrable random sequences $\mathbf{X}=\left(
X_{k}\right)  _{k\in\mathbb{Z}}$ with the first two conditional moments given
for all $k\in\mathbb{Z}$ by%

\begin{align}
\mathbb{E}\left(  X_{k}|...,X_{k-2},X_{k-1},X_{k+1},X_{k+2},...\right)   &
=L\left(  X_{k-1},X_{k+1}\right)  ,\label{linregr}\\
\mathbb{E}\left(  X_{k}^{2}|...,X_{k-2},X_{k-1},X_{k+1},X_{k+2},...\right)
&  =Q\left(  X_{k-1},X_{k+1}\right)  , \label{quadvar}%
\end{align}
where $L$ is a symmetric linear polynomial and $Q$ is a symmetric quadratic
polynomial. Such processes have been considered by Bryc \cite{bryc1}, who,
among other things, identified one-dimensional distributions of $\mathbf{X}$
for a special class of quadratic functions $Q$. In \cite{bryc2} Markov chains
with the same conditional structure has been investigated. In \cite{matszab} a
thorough analysis of the condition defining first conditional moments in
Bryc's paper \cite{bryc1} has been carried, which has made it possible to omit
some of the assumptions from \cite{bryc1} and simplify the proof of the result
contained therein.

The aim of this paper is to analyze random sequences satisfying (\ref{linregr}%
) and (\ref{quadvar}) without any restrictions on the parameters defining
polynomials $L$ and $Q$; the main result is Theorem \ref{mapa}, in which
description of almost all possible combinations of the parameters is given
and the distributions of the sequences are identified. One of the conclusions
of this analysis is partial negative answer to a question posed by Bryc in
\cite{bryc1}, whether there exist processes satisfying (\ref{linregr}) and
(\ref{quadvar}), whose one-dimensional distributions are not identified by
moments. The answer is based on some calculations involving Al-Salam-Chihara
polynomials, as it was shown in \cite{bms}. In addition, we reduce a number of
assumptions from \cite{bryc1}.

The paper is organized as follows. Subsection \ref{sectassumpt} lists the
assumptions used throughout the paper. Section \ref{sectmain} contains the
statements of the main results (their proofs are presented in section
\ref{sectproofs}) and some remarks. The next section contains auxiliary
results and their proofs.

\subsection{\label{sectassumpt}Assumptions and notation}

We will follow Bryc in using the term \textit{random field} for $\mathbf{X,}$
for additional comments see \cite{matszab}.

Let $\mathbf{X=}\left(  X_{k}\right)  _{k\in\mathbb{Z}}$ be a square
integrable random field indexed by the integers, with non-degenerate
covariance matrices and constant first two moments, that is $\mathbb{E}%
X_{k}=\mathbb{E}X_{0},$ $\mathbb{E}X_{k}^{2}=\mathbb{E}X_{0}^{2}$ $\forall
k\in\mathbb{Z}.$ After Bryc \cite{bryc1}, we assume that conditional structure
of $\mathbf{X}$ is given by (\ref{linregr}) and (\ref{quadvar}) with
\begin{align*}
L\left(  x,y\right)   &  =a\left(  x+y\right)  +b,\\
Q\left(  x,y\right)   &  =A\left(  x^{2}+y^{2}\right)  +Bxy+D\left(
x+y\right)  +C,
\end{align*}
$a,b,A,B,C,D\in\mathbb{R}.$ Non-singularity of covariance matrices implies
that all random variables $X_{k}$ are non-degenerate and there is no loss of
generality in assuming that $\mathbb{E}X_{k}=0$ and $\mathbb{E}X_{k}^{2}=1$
for all $k\in\mathbb{Z},$ which implies $b=0.$

It has been shown in \cite{matszab} that (\ref{linregr}) implies $L_{2}%
$-stationarity of $\mathbf{X.}$ Since the case $\rho:=\operatorname*{corr}%
\left(  X_{0},X_{1}\right)  =0$ contains sequences of independent random
variables (which satisfy (\ref{linregr}) and (\ref{quadvar}) but can have
arbitrary distributions), we shall exclude it from the considerations. Observe
that non-singularity of the covariance matrices implies $\left|  \rho\right|
<1.$ By Theorem 3.1 from \cite{bryc1} (see also Theorems 1 and 2 in
\cite{matszab}), $\operatorname*{corr}\left(  X_{0},X_{k}\right)
=\rho^{\left|  k\right|  }$ and one-sided regressions are linear%
\begin{equation}
\mathbb{E}\left(  X_{k}|...,X_{-1},X_{0}\right)  =\rho^{k}X_{0}.
\label{onesidedregr}%
\end{equation}
Multiplying (\ref{linregr}) by $X_{k-1}$ and taking the expected value, we see
that $a=\rho/\left(  1+\rho^{2}\right)  .$

Let us define $\sigma-$algebras $\mathcal{F}_{\leq m}:=\sigma\left(
X_{k}:k\leq m\right)  ,$ $\mathcal{F}_{\geq m}:=\sigma\left(  X_{k}:k\geq
m\right)  ,$ $\mathcal{F}_{\neq m}:=\sigma\left(  X_{k}:k\neq m\right)  .$
Throughout this paper $\operatorname*{supp}X$ stands for support of a random
variable $X$, $\#A$ denotes cardinality of a set $A,$ $\mathbb{I}_{A}$ is its
characteristic function, and $\delta_{p}$ is the Dirac measure concentrated on
$p\in\mathbb{R}.$

\section{\label{sectmain}Main results and an open problem}

In \cite{bryc1} it was assumed $D=0$ and
\begin{equation}
A\left(  \rho^{2}+\frac{1}{\rho^{2}}\right)  +B=1.\label{mysticformula}
\end{equation}

The following Theorem describes all valid combinations of parameters from
(\ref{quadvar}) (note that taking the expected value of (\ref{quadvar}) one
gets that $C=1-2A-B\rho^{2})$ and identifies the distributions of relevant fields.

\begin{theorem}\label{mapa}
There do not exist standardized random fields $\mathbf{X}=(X_{k})_{k\in\mathbb{Z}}$
with non-singular covariance matrices, satisfying (\ref{linregr}) and (\ref{quadvar})
and $\rho\neq0$, unless $D=0$ and either $B=0$ or (\ref{mysticformula}) holds with
$B\in{\mathcal{B}}_1\cup{\mathcal{B}}_2\cup{\mathcal{B}}_3$, where

\begin{gather}
{\mathcal{B}}_1=\left\{\frac{\rho^2-1}{\rho^2(1+\rho^2)}\right\},\
{\mathcal{B}}_2=\left(  0,\frac{2\rho^{2}}{\left(  1+\rho^{2}\right)  ^{2}}\right],\nonumber\\
{\mathcal{B}}_3=\left\{\frac{2\rho^{2}}{\left(  1+\rho^{2}\right)  ^{2}}\cdot\frac{\left(
1-\rho^{4}\right)  \left(  1+\rho^{2\slash m}\right)  }{2\left(  \rho^{2\slash
m}-\rho^{4}\right)}:m\in\mathbb{N}\right\}.\nonumber
\end{gather}
Furthermore, if $D=B=0$ or $D=0$ and (\ref{mysticformula}) holds with $B\in{\mathcal{B}}_2$ then
the relevant fields exist and have the following distributions:
\begin{enumerate}
\item  If $B=0$ and $A=1\slash2$ then $X_k^2=X_{k+1}^2$ for all $k$; if additionally there exists
$k_0\in\mathbf{Z}$ such that $\Pr(X_{k_0}=0)=0$ then $X_{k}=RY_{k},$ $R$ and $Y_{k}$ are
independent, $R\geq0,$ $\mathbb{E}R^{2}=1$ and all random variables $Y_{k}$ are symmetric with
values $\pm1;$

\item  If $B=0$ and $A\neq1\slash2$ then all random variables $X_{k}$ are
symmetric with values $\pm1;$

\item  If (\ref{mysticformula}) holds then for all $B\in\left(  0,2\rho
^{2}\slash\left(  1+\rho^{2}\right)  ^{2}\right)  $ all one-dimensional
distributions of $\mathbf{X}$ are equal to a uniquely determined symmetric
and absolutely continuous distribution (depending on $B$) supported on a finite interval;

\item  If (\ref{mysticformula}) holds and $B=2\rho^{2}\slash\left(1+\rho^{2}\right)^{2}$
then $X_{k}\sim\mathcal{N}\left(  0,1\right)  $ for all $k\in\mathbb{Z}$.

\end{enumerate}
\end{theorem}

\begin{theorem}\label{degenerat}
If (\ref{mysticformula}) holds with $D=0$ and $B\in{\mathcal{B}}_1$ then there
do not exist the relevant random fields $\mathbf{X}$ satisfying
(\ref{linregr}-\ref{quadvar}) with $\rho\neq0$ and such that
sequence $(X^2_k)_k$ is uniformly integrable and $\Pr(X_{k-1}=0)<1$.
\end{theorem}

\begin{problem}
 Do there exist random fields $\mathbf{X}$ satisfying
(\ref{linregr}-\ref{quadvar}) with $\rho\neq0$, $D=0$, and such that (\ref{mysticformula}) holds
with $B\in{\mathcal{B}}_3$?
\end{problem}

Theorem \ref{markov} strengthens Theorem 2.3 from \cite{bryc2} by stating that
$\mathbf{X}$ is always a Markov chain and relaxing assumptions for the
uniquely determined case.

\begin{theorem}\label{markov}
If $\mathbf{X}$ is a standardized random field with non-singular
covariance matrices, satisfying (\ref{linregr}) and (\ref{quadvar}), $\rho\neq0$, and such that
$D=B=0$ or $D=0$ and (\ref{mysticformula}) holds with $B\in{\mathcal{B}}_2$,
then $\mathbf{X}$ is a stationary Markov chain. $\mathbf{X}$ has
uniquely determined finite-dimensional distributions unless $B=D=0$ and $A=1/2$.
\end{theorem}

Clearly, if $A=1\slash2,$ $B=D=0$ then distributions of $\mathbf{X}$ can be
arbitrary (but symmetric). Observe that in this case condition (\ref{quadvar})
is trivial.

\subsection{Remarks}

\begin{enumerate}
\item  Unlike Bryc \cite{bryc1}, we do not assume $L_{2}$-stationarity of
$\mathbf{X}$ (which is proved in \cite{matszab}), nor do we assume equality of
its one-dimensional distributions (it is stated in Theorem \ref{mapa}).
However, to give a simple and uniform description of the considered random
fields, we shall apply some results of Bryc from \cite{bryc1}. One can easily
verify that the results from \cite{bryc1}, used in sections \ref{sectauxi} and
\ref{sectproofs}, do only use assumptions listed in subsection
\ref{sectassumpt} (in particular, they do not rely on the equality of
one-dimensional distributions of $\mathbf{X}$).

\item  We do not use assumption (6) from \cite{bryc1}. The only restriction on
the correlation coefficients (apart from those stemming from non-singularity
of covariance matrices - see \cite{matszab}, Remark 2) will be requirement
that $\rho\neq0.$

\item Bryc \cite{bryc1} assumed $D=0$ and (\ref{mysticformula}) and introduced
new parameters%
\begin{equation}
R:=B{\left(\rho+\frac{1}{\rho}\right)}^2  ,\quad q:=\frac{\rho^{4}%
+R-1}{1+\rho^{4}\left(  R-1\right)  }.\label{q}%
\end{equation}
Using some facts concerning $q$-Hermite polynomials he proved (under
additional assumption of equality of one--dimensional distributions, see
\cite{bryc1}, Theorem 3.2) implications (3) and (4) of the above Theorem and
the fact that then there do not exist $\mathbf{X}$ with $R<0$ (equivalently
$B<0$). (If $R=B=0$ then Bryc showed that one--dimensional distributions are
symmetric with values $\pm1,$ what can be deduced from Theorem \ref{mapa},
since then (\ref{mysticformula}) implies $A=\rho^{2}\slash\left(  1+\rho
^{4}\right)  \neq1\slash2$.) Bryc also raised a question (Concluding remarks,
4) whether there exist stationary processes which satisfy conditions
(\ref{linregr}), (\ref{quadvar}), (\ref{mysticformula}) and $D=0,$
corresponding to $R>2$ (equivalently $B>2\rho^{2}\slash\left(  1+\rho
^{2}\right)  ^{2}$). As Theorem \ref{mapa} shows, the answer is almost
negative. More precisely there exists a discrete set of possible values of $B$
all contained in the interval $(2\rho^{2}\slash\left(  1+\rho^{2}\right)
^{2},1]$ for which one--dimensional distributions of $\mathbf{X}$ have
all moments but are not determined by moments and all conditional
distributions $X_{j}|X_{k}$, $j\neq k$ are discrete.

\item  The densities from implication (3) in Theorem \ref{mapa} have explicit
product representation - see \cite{bryc1}, Concluding remarks 3.

\item Observe that the case considered in Theorem \ref{degenerat} leads to undefined value
of $q$ (see (\ref{q})). Therefore it has to be examined separately. In this case (unlike in cases
in which $q$ as a function of $B$ is determined), we cannot apply Corrolary 5.1 from \cite{bryc1}
and deduce existence of higher moments of $\mathbf{X}$. If we knew that all moments exist, then
of course we could have dropped the uniform integrability assumption.

\end{enumerate}

\section{\label{sectauxi}Auxiliary results}

\begin{proposition}
\label{twopoint}Let $A\neq1\slash\left(  1+\rho^{2}\right)  .$ If
$\#\operatorname*{supp}X_{k}=2$ for all $k$ then $D=B=0$ and $X_{k}\sim\left(
\delta_{-1}+\delta_{1}\right)  \slash2,$ or (\ref{mysticformula}) holds and $D=0.$
\end{proposition}

\begin{proof}
From the proof of Lemma 5.2, \cite{bryc1} (note that there is a minor misprint
in \cite{bryc1} in the formula for $\mathbb{E}\left(  X_{1}^{2}|X_{0}\right)
$, repeated in the proof of Proposition 3.1) it follows that
\begin{align}
\mathbb{E}\left(  X_{k+1}^{2}|\mathcal{F}_{\leq k}\right)   &  =\alpha
_{1}X_{k}^{2}+\beta_{1}X_{k}+\gamma_{1},\label{w1k1}\\
\mathbb{E}\left(  X_{k+2}^{2}|\mathcal{F}_{\leq k}\right)   &  =\alpha
_{2}X_{k}^{2}+\beta_{2}X_{k}+\gamma_{2}\label{w1k2}%
\end{align}
with%
\begin{align*}
\alpha_{1} &  =\frac{A\left(  1-\rho^{2}\right)  +B\rho^{2}}{1-A\left(
1+\rho^{2}\right)  },\alpha_{2}=\frac{\left(  1+\rho^{2}\right)  \left(
A+B\rho^{2}\right)  -\rho^{2}}{1-A\left(  1+\rho^{2}\right)  },\\
\beta_{1} &  =\frac{D\left(  1+\rho^{2}\right)  }{1-A\left(  1+\rho
^{2}\right)  },\beta_{2}=\frac{D\left(  1+\rho^{2}\right)  ^{2}}{1-A\left(
1+\rho^{2}\right)  },\\
\gamma_{1} &  =\frac{C}{1-A\left(  1+\rho^{2}\right)  },\gamma_{2}%
=\frac{C\left(  1+\rho^{2}\right)  }{1-A\left(  1+\rho^{2}\right)  }.
\end{align*}
Observe that $\mathbb{E}\left(  X_{k+2}^{2}|\mathcal{F}_{\leq k}\right)
=\mathbb{E}\left[  \mathbb{E}\left(  X_{k+2}^{2}|\mathcal{F}_{\leq
k+1}\right)  |\mathcal{F}_{\leq k}\right]  ,$ therefore for all $k\in
\mathbb{Z}$
\begin{equation}
\alpha_{2}X_{k}^{2}+\beta_{2}X_{k}+\gamma_{2}=\alpha_{1}^{2}X_{k}^{2}%
+X_{k}\beta_{1}\left(  \alpha_{1}+\rho\right)  +\gamma_{1}\left(  \alpha
_{1}+1\right)  .\label{rkwadratowe}%
\end{equation}

First assume that
\begin{equation}
\left\{
\begin{array}
[c]{c}%
\alpha_{1}^{2}=\alpha_{2}\\
\beta_{1}\left(  \alpha_{1}+\rho\right)  =\beta_{2}\\
\gamma_{1}\left(  \alpha_{1}+1\right)  =\gamma_{2}%
\end{array}
\right.  .\label{rkwadratowe1}%
\end{equation}
This after some easy algebra gives
\begin{align}
D\text{\thinspace}\left\{  C-\left(  1-\rho\right)  \,\rho\,\left[
1-A\,\left(  1+{{\rho}^{2}}\right)  \right]  \right\}   &  =0,\label{wzor1}\\
C{\,\left[  A\left(  \rho^{2}+\frac{1}{\rho^{2}}\right)  +B\,{-1}\right]  } &
={0.}\label{wzor2}%
\end{align}

If $C{=0}$ {then }$D=0$ by (\ref{wzor1}) and $A=\left(  1-B\rho^{2}\right)
/2$. From (\ref{w1k1}), $\mathbb{E}\left(  X_{k+1}^{2}|\mathcal{F}_{\leq
k}\right)  =X_{k}^{2}.$ Consequently, by symmetry, $\mathbb{E}\left(
X_{k}^{2}|\mathcal{F}_{\geq k+1}\right)  =X_{k+1}^{2},$ hence $\mathbb{E}%
\left(  X_{k+1}^{2}|X_{k}^{2}\right)  =X_{k}^{2}$ and $\mathbb{E}\left(
X_{k}^{2}|X_{k+1}^{2}\right)  =X_{k+1}^{2}.$ By a result of Doob (see
\cite{doob}, page 314) $X_{k}^{2}=X_{k+1}^{2},$ so (since $\#\operatorname*{supp}X_{k}=2$) $X_{k}^{2}=1$ and
$X_{k}\sim\left(  \delta_{-1}+\delta_{1}\right)  \slash2$ for all $k.$ From
(\ref{quadvar}) $B\left(  X_{k-1}X_{k+1}-\rho^{2}\right)  =0,$ so $B=0$ by the
assumptions on covariance matrices and $A=1\slash2.$

If $C\neq0$ then (\ref{wzor2}) implies (\ref{mysticformula}) and (\ref{wzor1})
takes the form
\[
D\left(  1-\rho\right)  \left(  1+\rho^{2}\right)  \left(  \frac{1}{1+\rho
^{2}}-A\right)  =0,
\]
hence $D=0.$

Now assume that (\ref{rkwadratowe1}) does not hold. Then the values of $X_{k}$
have to be the roots of (\ref{rkwadratowe}), so all $X_{k}$ have to be
equidistributed. Observe that $\mathbb{E}X_{k}^{3}=0$ for all $k\in
\mathbb{Z}.$ Indeed, applying (\ref{onesidedregr}) twice one gets
\begin{multline*}
\mathbb{E}\left(  X_{k-1}X_{k}X_{k+1}\right)  =\mathbb{E}\left[  X_{k-1}%
X_{k}\mathbb{E}\left(  X_{k+1}|\mathcal{F}_{\leq k}\right)  \right]
=\rho\mathbb{E}\left(  X_{k-1}X_{k}^{2}\right)  \\
=\rho\mathbb{E}\left[  \mathbb{E}\left(  X_{k-1}|\mathcal{F}_{\geq k}\right)
X_{k}^{2}\right]  =\rho^{2}\mathbb{E}X_{k}^{3}%
\end{multline*}
and on the other hand, by using (\ref{linregr}) and (\ref{onesidedregr})
\[
\mathbb{E}\left(  X_{k-1}X_{k}X_{k+1}\right)  =\mathbb{E}\left[
X_{k-1}\mathbb{E}\left(  X_{k}|\mathcal{F}_{\neq k}\right)  X_{k+1}\right]
=\frac{\rho^{3}}{1+\rho^{2}}\left(  \mathbb{E}X_{k-1}^{3}+\mathbb{E}%
X_{k+1}^{3}\right)  .
\]
Thus $X_{k}\sim\left(  \delta_{-1}+\delta_{1}\right)  \slash2$ and from
(\ref{quadvar})
\begin{equation}
B\left(  X_{k-1}X_{k+1}-\rho^{2}\right)  +D\left(  X_{k-1}+X_{k+1}\right)
=0.\label{allberno1}%
\end{equation}
Define $F_{1}:=\left\{  X_{k-1}=-X_{k+1}\right\}  ,$ $F_{2}:=\left\{
X_{k-1}=X_{k+1}\right\}  .$ By the assumptions on covariance matrices
$\Pr\left(  F_{1}\right)  \neq1$ and $\Pr\left(  F_{2}\right)  \neq1.$
Applying (\ref{allberno1}) to $\omega\in F_{1}$ and $\omega\in F_{2},$ we get
$B=0$ and $D=0$. Note that (\ref{quadvar}) reduces to a trivial identity.
\end{proof}

\begin{proposition}
\label{manypoint}Let $A\neq1\slash\left(  1+\rho^{2}\right)  .$ If there exists $k_0\in\mathbb{Z}$
such that $\#\operatorname*{supp}X_{k_0}>2$ then either $A=1\slash2,B=C=D=0$ and $X_k^2=X_{k+1}^2$
for all $k$, or $(\ref{mysticformula})$ holds and $D=0\neq C.$ If, additionally, there exists
$k_1\in\mathbb{Z}$ such that $\Pr(X_{k_1}=0)=0$ then in the first case we have $X_{k}=RY_{k},$
where $R\geq0,$ $Y_{k}\sim\left( \delta_{-1}+\delta_{1}\right)  /2$ for all $k$ and $R$ and $Y_{k}$
are independent.
\end{proposition}

\begin{proof}
We keep the notation from the proof of Proposition \ref{twopoint}. If there
exists $k_0\in\mathbb{Z}$ such that $\#\operatorname*{supp}X_{k_0}>2$ then
(\ref{rkwadratowe1}) must be satisfied, so (\ref{wzor1}) and (\ref{wzor2}) hold.

If $C{=0}$ {then }$X_{k}^{2}=X_{k+1}^{2}$ for all $k\in\mathbb{Z}$ as in the proof of Proposition
\ref{twopoint}. By (\ref{quadvar}),
\begin{equation}
X_{k-1}\left[  \left(  2A-1\right)  X_{k-1}+BX_{k+1}\right]
=0.\label{wzorek1}%
\end{equation}
Observe that $\left\{  X_{k-1}=0\right\}  =\left\{  X_{k+1}=0\right\}  ,$ so
if $A\neq1/2$ (hence $B\neq0)$ then $X_{k+1}=X_{k-1}\left(  1-2A\right)  /B,$
which contradicts the non-singularity of the covariance matrices. Note that if
$A=1/2$ then (\ref{quadvar}) is automatically fulfilled. Since $X_{k}%
^{2}=X_{k+1}^{2}=:R^{2},$ where $R\geq0,$ then, under the assumption of lack of atom at $0$,
$X_{k}=RY_{k}$, where $\operatorname*{supp}Y_{k}=\left\{  -1,1\right\}$. By (\ref{onesidedregr}),
$\mathbb{E}\left(  X_{k}|X_{k-1}\right)  =R\ \mathbb{E}\left(  Y_{k}%
|X_{k-1}\right)  =\rho R\ Y_{k-1},$ so $\mathbb{E}\left(  Y_{k}|X_{k-1}%
\right)  =\rho Y_{k-1}$ and $\mathbb{E}\left(  Y_{k}|R\right)  =\rho
\mathbb{E}\left(  Y_{k-1}|R\right)  .$ Analogously, $\mathbb{E}\left(
Y_{k-1}|X_{k}\right)  =\rho Y_{k}$ and $\mathbb{E}\left(  Y_{k-1}|R\right)
=\rho\mathbb{E}\left(  Y_{k}|R\right)  ,$ therefore $\mathbb{E}\left(
Y_{k-1}|R\right)  =\mathbb{E}\left(  Y_{k}|R\right)  =0.$ Hence $\mathbb{E}%
\left(  Y_{k}\right)  =0$ and $\Pr\left(  Y_{k}=1|R\right)  =\Pr\left(
Y_{k}=-1|R\right)  =1/2.$ Thus $Y_{k}$ and $R$ are independent and $Y_{k}%
\sim\left(  \delta_{-1}+\delta_{1}\right)  /2$.

If $C\neq0$ then (\ref{wzor2}) implies (\ref{mysticformula}) and $D=0$ as in
the proof of Proposition \ref{twopoint}.
\end{proof}

The following Lemma  will be used in the proof of Theorem \ref{degenerat}.

\begin{lemma}
Suppose that (\ref{mysticformula}) holds with $D=0$ and $B\in{\mathcal{B}}_1$.
If $\mathcal{G}_{n}:=\sigma\left(  \mathcal{F}_{\leq k-1},\mathcal{F}_{\geq
k+n}\right)  $ then for all $n=1,2,...$ and $k\in\mathbb{Z}$%
\begin{equation}
\mathbb{E}\left(  X_{k}^{2}|\mathcal{G}_{n}\right)  =\frac{1-\rho^{2n}}%
{1-\rho^{2n+2}}X_{k-1}^{2}+\frac{1-\rho^{2}}{1-\rho^{2n+2}}X_{k+n}^{2}%
+\frac{\left(  \rho^{2}-1\right)  \left(  1-\rho^{2n}\right)  }{\rho
^{n+1}\left(  1-\rho^{2n+2}\right)  }X_{k-1}X_{k+n}.\label{wzor5}%
\end{equation}
\end{lemma}

\begin{proof}
The proof is by induction on $n.$ (\ref{wzor5}) for $n=1$ is (\ref{quadvar}).

Assume that the assertion holds for some $n\geq1.$ Then
\begin{multline*}
\mathbb{E}\left(  X_{k}^{2}|\mathcal{G}_{n+1}\right)  =\mathbb{E}\left[
\mathbb{E}\left(  X_{k}^{2}|\mathcal{F}_{\neq k}\right)  |\mathcal{G}%
_{n+1}\right]  =\\
=AX_{k-1}^{2}+A\mathbb{E}\left(  X_{k+1}^{2}|\mathcal{G}_{n+1}\right)
+BX_{k-1}\mathbb{E}\left(  X_{k+1}|\mathcal{G}_{n+1}\right)  .
\end{multline*}
By the induction assumption
\begin{multline*}
\mathbb{E}\left(  X_{k+1}^{2}|\mathcal{G}_{n+1}\right)  =\mathbb{E}\left[
\mathbb{E}\left(  X_{k+1}^{2}|\sigma\left(  \mathcal{F}_{\leq k}%
,\mathcal{F}_{\geq k+n+1}\right)  \right)  |\mathcal{G}_{n+1}\right]  =\\
=\frac{1-\rho^{2n}}{1-\rho^{2n+2}}\mathbb{E}\left(  X_{k}^{2}|\mathcal{G}%
_{n+1}\right)  +\frac{1-\rho^{2}}{1-\rho^{2n+2}}X_{k+n+1}^{2}+\\
+\frac{\left(  \rho^{2}-1\right)  \left(  1-\rho^{2n}\right)  }{\rho
^{n+1}\left(  1-\rho^{2n+2}\right)  }X_{k+n+1}\mathbb{E}\left(  X_{k}%
|\mathcal{G}_{n+1}\right)  .
\end{multline*}
This, formulae
\begin{align*}
\mathbb{E}\left(  X_{k}|\mathcal{G}_{n+1}\right)   &  =\frac{\rho^{n+1}%
-\rho^{-n-1}}{\rho^{n+2}-\rho^{-2-n}}X_{k-1}+\frac{\rho-\frac{1}{\rho}}%
{\rho^{n+2}-\rho^{-2-n}}X_{k+n+1},\\
\mathbb{E}\left(  X_{k+1}|\mathcal{G}_{n+1}\right)   &  =\frac{\rho^{n}%
-\rho^{-n}}{\rho^{n+2}-\rho^{-2-n}}X_{k-1}+\frac{\rho^{2}-\frac{1}{\rho^{2}}%
}{\rho^{n+2}-\rho^{-2-n}}X_{k+n+1}%
\end{align*}
(see \cite{matszab}, formula (3.4) in the proof of Theorem 2) and
an easy observation that $1-A\left(  1-\rho^{2n}\right)  /\left(
1-\rho^{2n+2}\right)  \neq0$ for all $n\in\mathbb{N},$ completes the proof.
\end{proof}

\section{\label{sectproofs}Proofs of main results}

\begin{proof}
[Proof of Theorem \ref{mapa}]A starting point of our considerations is Lemma
5.2, formula (16), \cite{bryc1}:
\begin{equation}
\mathbb{E}\left(  X_{k+1}^{2}|\mathcal{F}_{\leq k}\right)  \left[  1-A\left(
1+\rho^{2}\right)  \right]  =X_{k}^{2}\left[  A\left(  1-\rho^{2}\right)
+B\rho^{2}\right]  +X_{k}D\left(  1+\rho^{2}\right)  +C\label{bryclemma52}%
\end{equation}
(in fact, the above formula in \cite{bryc1} is \ given for $k=0$ only but from
its proof it is evident that (\ref{bryclemma52}) is true for all
$k\in\mathbb{Z}$).

If $A=1\slash\left(  1+\rho^{2}\right)  $ then (\ref{bryclemma52}) is a
quadratic equation and if $A\left(  1-\rho^{2}\right)  +B\rho^{2}\neq0$,
equivalently $B\neq\left(  \rho^{2}-1\right)  \slash\left[  \rho^{2}\left(
\rho^{2}+1\right)  \right]  ,$ it has two solutions, therefore all $X_{k}$
must have the same two point distribution. Repeating the reasoning from the
proof of Proposition \ref{twopoint}, we see that $\mathbb{E}X_{k}^{3}=0$ for
all $k,$ so $X_{k}\sim\left(  \delta_{-1}+\delta_{1}\right)  /2$ and $B=D=0.$

If $A=1/\left(  1+\rho^{2}\right)  $ and $B=\left(  \rho^{2}-1\right)
\slash\left[  \rho^{2}\left(  \rho^{2}+1\right)  \right]  $ then
$C=1-2A-B\rho^{2}=0$ and $D=0$ by (\ref{bryclemma52}). This case is analyzed in Theorem
\ref{degenerat}.

Now assume $A\neq1\slash\left(  1+\rho^{2}\right)  .$ By Propositions
\ref{twopoint} and \ref{manypoint}, it suffices to examine the case when
(\ref{mysticformula}) and $D=0$ hold. (This is the combination of parameters
assumed by Bryc in \cite{bryc1}.) In this case, $\mathbb{E}\left|
X_{k}\right|  ^{p}<\infty$ for all $p>1$ by \cite{bryc1}, Corollary 5.1. By
the proof of \cite{bryc1}, Lemma 6.3, for all $k\in\mathbb{Z}$ and $n\geq1$%
\begin{equation}
\mathbb{E}\left[  Q_{n}\left(  X_{k}\right)  |\mathcal{F}_{\leq k-1}\right]
=\rho^{n}Q_{n}\left(  X_{k-1}\right)  ,\label{lemma63}%
\end{equation}
where $Q_{n}$ are $q$-Hermite polynomials defined in \cite{bryc1}, Subsection
6.1, with $q$ as in (\ref{q}). (\ref{lemma63}) and symmetry imply
\[
\mathbb{E}\left[  Q_{n}\left(  X_{k}\right)  |X_{k-1}\right]  =\rho^{n}%
Q_{n}\left(  X_{k-1}\right)  ,\quad\mathbb{E}\left[  Q_{n}\left(
X_{k-1}\right)  |X_{k}\right]  =\rho^{n}Q_{n}\left(  X_{k}\right)  ,
\]
hence for $n,m\in\mathbb{N}$ and $k\in\mathbb{Z}$ we have
\begin{align*}
\mathbb{E}\left[  Q_{n}\left(  X_{k}\right)  Q_{m}\left(  X_{k-1}\right)
\right]   &  =\rho^{m}\mathbb{E}\left[  Q_{n}\left(  X_{k}\right)
Q_{m}\left(  X_{k}\right)  \right]  =\rho^{n}\mathbb{E}\left[  Q_{n}\left(
X_{k-1}\right)  Q_{m}\left(  X_{k-1}\right)  \right]  ,\\
\mathbb{E}\left[  Q_{n}\left(  X_{k-1}\right)  Q_{m}\left(  X_{k}\right)
\right]   &  =\rho^{m}\mathbb{E}\left[  Q_{n}\left(  X_{k-1}\right)
Q_{m}\left(  X_{k-1}\right)  \right]  =\rho^{n}\mathbb{E}\left[  Q_{n}\left(
X_{k}\right)  Q_{m}\left(  X_{k}\right)  \right]  .
\end{align*}
Denoting $x:=\mathbb{E}\left[  Q_{n}\left(  X_{k}\right)  Q_{m}\left(
X_{k}\right)  \right]  ,$ $y:=\mathbb{E}\left[  Q_{n}\left(  X_{k-1}\right)
Q_{m}\left(  X_{k-1}\right)  \right]  $ one obtains $\rho^{m}x=\rho^{n}y$ and
$\rho^{m}y=\rho^{n}x,$ hence (if $\ n\neq m$)$,$ $x=y=0$ and $Q_{n}$ are
orthogonal with respect to $\mathcal{L}\left(  X_{k-1}\right)  $ and
$\mathcal{L}\left(  X_{k}\right)  .$ By Lemma 6.1, \cite{bryc1} we deduce that
there do not exist random fields with (\ref{mysticformula}), $D=0$ and $B<0$
(equivalently $q<-1$) and (using also Proposition 8.1 from \cite{bryc1}) that
assertions 3. and 4. of Theorem \ref{mapa} are true (by
Lemma 6.1, \cite{bryc1}, in cases covered by assertions 3. and 4.,
measures that make $Q_{n}$ orthogonal are unique, so one-dimensional
distributions of $\mathbf{X}$ are equal).

What is left is to exclude the case of (\ref{mysticformula}), $D=0$ and
$B>2\rho^{2}\slash\left(  1+\rho^{2}\right)  ^{2}$ (equivalently $q>1$). (This
case was unresolved in \cite{bryc1}.) To obtain a contradiction, suppose that
there exists $\mathbf{X}$ with this set of parameters and consider the
conditional distribution $\mathcal{L}\left(  X_{1}|X_{0}=y\right)  .$ By
\cite{bms}, Theorem 2, its monic orthogonal polynomials are Al-Salam--Chihara
polynomials (see \cite{bms} for more details) and satisfy the three term
recurrence relation%
\[
p_{n+1}\left(  x\right)  =\left(  x-\rho yq^{n}\right)  p_{n}\left(  x\right)
-\left(  1-\rho^{2}q^{n-1}\right)  \left[  n\right]  _{q}p_{n-1}\left(
x\right)  ,
\]
where $\left[  n\right]  _{q}=1+q+...+q^{n-1}.$ Since the distribution
$\mathcal{L}\left(  X_{1}|X_{0}=y\right)  $ is a positive measure, the
coefficients at the third term in the above recurrence must be non-negative, thus
$(1-\rho^{2}q^{n-1})\geq0$ for all $n\in\mathbb{N},$ so if $q\ne\rho^{-2/m}$ for all $m\in\mathbb{N}$
we obtain a contradiction.

If $q=\rho^{-2/m}$ for some $m\in\mathbb{N}$ then conditional
$\mathcal{L}\left(  X_{1}|X_{0}=y\right)  $ and two--dimensional
$\mathcal{L}\left(  X_{0},X_{1}\right)$ distributions do exist. The question whether there exist relevant
random fields remains open.
\end{proof}

\begin{proof}[Proof of Theorem \ref{degenerat}]
Fix $k\in\mathbb{Z}$. Since by L\'{e}vy's Downward Theorem sequence $\left\{  \mathbb{E}\left(
X_{k}^{2}|\mathcal{G}_{n}\right)  \right\}  _{n=1}^{+\infty}$ converges almost
surely and in $L_{1},$ multiplying both sides of (\ref{wzor5}) by the
indicator of the event $\left\{  X_{k-1}\neq0\right\}  $ and $\rho^{n+1}$ and
letting $n\rightarrow+\infty$ we see that $X_{n}\rightarrow0$ for almost all
$\omega.$ \ Now define $Y_{n}:=\mathbb{E}\left(  X_{n+k}^{2}|\mathcal{F}_{\leq
k}\right)  .$ Then%
\begin{align*}
Y_{n}  & =\mathbb{E}\left(  \mathbb{E}\left(  X_{n+k}^{2}|\mathcal{F}_{\neq
n+k}\right)  |\mathcal{F}_{\leq k}\right)  \\
& =\frac{1}{1+\rho^{2}}\left(  \mathbb{E}\left(  X_{k+n-1}^{2}|\mathcal{F}%
_{\leq k}\right)  +\mathbb{E}\left(  X_{k+n+1}^{2}|\mathcal{F}_{\leq
k}\right)  \right)  -\frac{1-\rho^{2}}{\rho^{2}\left(  1+\rho^{2}\right)
}\mathbb{E}\left(  X_{k+n-1}X_{k+n+1}|\mathcal{F}_{\leq k}\right)  \\
& =\left(  \frac{1}{1+\rho^{2}}-\frac{1-\rho^{2}}{\left(  1+\rho^{2}\right)
}\right)  \mathbb{E}\left(  X_{k+n-1}^{2}|\mathcal{F}_{\leq k}\right)
+\frac{1}{1+\rho^{2}}\mathbb{E}\left(  X_{k+n+1}^{2}|\mathcal{F}_{\leq
k}\right)  \\
& =\frac{\rho^{2}}{1+\rho^{2}}\mathbb{E}\left(  X_{k+n-1}^{2}|\mathcal{F}%
_{\leq k}\right)  +\frac{1}{1+\rho^{2}}\mathbb{E}\left(  X_{k+n+1}%
^{2}|\mathcal{F}_{\leq k}\right)  =\frac{\rho^{2}}{1+\rho^{2}}Y_{n-1}+\frac
{1}{1+\rho^{2}}Y_{n+1}.
\end{align*}
Therefore we have a recurrence with obvious initial conditions $Y_{0}%
=X_{k}^{2}$ and $Y_{1}=\mathbb{E}\left(  X_{k+1}^{2}|\mathcal{F}_{\leq
k}\right)  .$ Since its characteristic equation is
\[
x^{2}-\left(  1+\rho^{2}\right)  x+\rho^{2}=0,
\]
we have
\[
Y_{n}=C_{1}+C_{2}\rho^{2n}%
\]
for some $\mathcal{F}_{\leq k}$-measurable random variables $C_{1}$ and
$C_{2}.$ Taking the initial conditions into account we obtain%
\[
C_{2}=\frac{X_{k}^{2}-\mathbb{E}\left(  X_{k+1}^{2}|\mathcal{F}_{\leq
k}\right)  }{1-\rho^{2}},C_{1}=-\frac{\rho^{2}X_{k}^{2}-\mathbb{E}\left(
X_{k+1}^{2}|\mathcal{F}_{\leq k}\right)  )}{1-\rho^{2}},
\]
so
\[
Y_{n}=X_{k}^{2}+\frac{X_{k}^{2}-\mathbb{E}\left(  X_{k+1}^{2}|\mathcal{F}%
_{\leq k}\right)  }{1-\rho^{2}}\left(  \rho^{2n}-1\right)  .
\]
From the uniform integrability assumption and the fact that $X_{n}$ converges
almost surely to $0,$
\[
X_{k}^{2}-\frac{X_{k}^{2}-\mathbb{E}\left(  X_{k+1}^{2}|\mathcal{F}_{\leq
k}\right)  }{1-\rho^{2}}=0,
\]
thus $\mathbb{E}\left(  X_{k+1}^{2}|\mathcal{F}_{\leq k}\right)  =\rho
^{2}X_{k}^{2}$. \ Since it implies $\operatorname*{var}\left(  X_{k+1}|\mathcal{F}_{\leq k}\right)  =0,$
we have a contradiction,
\end{proof}

\begin{proof}
[Proof of Theorem \ref{markov}]By Theorem \ref{mapa}, we have to examine three cases:

\begin{enumerate}
\item (\ref{mysticformula}) holds with $D=0$ and $C\neq0,$

\item $A\neq1/2,$ $B=D=0,C=1-2A\neq0,$

\item $A=1/2,$ $B=C=D=0.$
\end{enumerate}

In the first case $B\in\left(  0,2\rho^{2}/\left(  1+\rho^{2}\right)
^{2}\right]  ,$ what is equivalent to $\rho^{2}/\left(  1+\rho^{2}\right)
^{2}\leq A<\rho^{2}/\left(  1+\rho^{4}\right)  ,$ so in this case Markov
property was proved in Theorem 2.3, \cite{bryc2}. We shall present a different
approach, based on (\ref{lemma63}). Of course, it suffices to prove that for all $k\in\mathbb{Z}$ and
$B\in\mathcal{B}\left(  \mathbb{R}\right)  ,$ $\Pr\left(  X_{k+1}\in
B|\mathcal{F}_{\leq k}\right)  =\Pr\left(  X_{k+1}\in B|\mathcal{F}%
_{=k}\right)  .$ Fix $k\in\mathbb{Z}$, $B\in\mathcal{B}\left(  \mathbb{R}%
\right)  $ and define $\alpha_{j}:=\mathbb{E}\mathbb{I}_{B}\left(
X_{k+1}\right)  Q_{j}\left(  X_{k+1}\right)  ,$ where polynomials $Q_{j}$ are
$q$-Hermite polynomials with $q$ defined by (\ref{q}). M. Riesz proved that if
a positive measure $\mu$ on $\mathbb{R}$ is determined by moments then
$\mathbb{C}\left[  x\right]  $ is dense in $L_{2}\left(  \mu\right)  $ (see
\cite{akhiezer}, Corollary 2.3.3). This implies
\begin{equation}
\mathbb{I}_{B}\left(  X_{k+1}\right)  \overset{L_{2}}{=}\sum_{j=0}^{\infty
}\alpha_{j}Q_{j}\left(  X_{k+1}\right)  .\label{indyk}%
\end{equation}
Using (\ref{lemma63}) one obtains that%
\begin{multline*}
\mathbb{E}\left|  \mathbb{I}_{B}\left(  X_{k+1}\right)  -\sum_{j=0}^{n}%
\alpha_{j}Q_{j}\left(  X_{k+1}\right)  \right|  ^{2}\geq\\
\geq\mathbb{E}\left|  \mathbb{E}\left[  \left.  \mathbb{I}_{B}\left(
X_{k+1}\right)  -\sum_{j=0}^{n}\alpha_{j}Q_{j}\left(  X_{k+1}\right)  \right|
\mathcal{F}_{\leq k}\right]  \right|  ^{2}=
\end{multline*}%
\[
=\mathbb{E}\left|  \mathbb{E}\left(  \mathbb{I}_{B}\left(  X_{k+1}\right)
|\mathcal{F}_{\leq k}\right)  -\sum_{j=0}^{n}\alpha_{j}\rho^{j}Q_{j}\left(
X_{k}\right)  \right|  ^{2}.
\]
Passing to the limit as $n\rightarrow\infty$ and using (\ref{indyk}) we obtain
that $\mathbb{E}\left(  \mathbb{I}_{B}\left(  X_{k+1}\right)  |\mathcal{F}%
_{\leq k}\right)  =\sum_{j=0}^{\infty}\alpha_{j}\rho^{j}Q_{j}\left(
X_{k}\right)  $ in $L_{2}.$ If $\sum\alpha_{j}^{2}<+\infty$ then $\sum
\alpha_{j}^{2}\rho^{2j}\log^{2}j<+\infty,$ hence
\[
\Pr\left(  X_{k+1}\in B|\mathcal{F}_{\leq k}\right)  =\mathbb{E}\left(
\mathbb{I}_{B}\left(  X_{k+1}\right)  |\mathcal{F}_{\leq k}\right)
=\sum_{j=0}^{\infty}\alpha_{j}\rho^{j}Q_{j}\left(  X_{k}\right)
\]
almost surely by the Rademacher-Menshov theorem. Thus $\Pr\left(  X_{k+1}\in
B|\mathcal{F}_{\leq k}\right)  =\Pr\left(  X_{k+1}\in B|\mathcal{F}%
_{=k}\right)$ almost surely.

Now consider $A\neq1/2,$ $B=D=0,$ $C=1-2A$. In this case all random variables
are symmetric with values $\pm1,$ so this case is the same as
(\ref{mysticformula}), $D=0$, $A=\rho^{2}/\left(  1+\rho^{4}\right)  ,$ and
Markov property follows from Theorem 2.3, \cite{bryc2}.

Now consider $A=1/2,$ $B=C=D=0.$ Then $X_{k}=RY_{k},$ with $R\geq0$ and all
random variables $Y_{k}$ symmetric with values $\pm1.$ Since $\Pr\left(
R\geq0\right)  =1$ and by (\ref{onesidedregr}), $\mathbb{E}\left(
X_{k}|X_{k-1}\right)  =R\mathbb{E}\left(  Y_{k}|X_{k-1}\right)  =\rho
RY_{k-1},$ so $\mathbb{E}\left(  Y_{k}|Y_{k-1}\right)  =\rho Y_{k-1},$ hence
$\left(  Y_{k}\right)  _{k}$ is a Markov chain, again by Bryc's argument.
Since $R$ and $Y_{k}$ are independent and $\left(  R,Y\right)  \mapsto RY$ is
one-to-one, $\mathbf{X}$ is a Markov chain.

In the first two cases all finite-dimensional distributions are uniquely
determined by Theorem 2.3, \cite{bryc2}. Clearly, one cannot determine
distributions of $\mathbf{X}$ when the third case holds.
\end{proof}

\textbf{Acknowledgement }\textit{The authors thank W. Bryc, R. Lata\l a and J. Weso\l owski
for several very helpful discussions.}


\begin{thebibliography}{99}
\bibitem{akhiezer}Akhiezer, N.I. (1965) The Classical Moment Problem. Oliver
and Boyd, Edinburgh. Translated from the Russian.

\bibitem {bryc1}Bryc, W. (2001) Stationary random fields with linear
regressions. Annals of Probability 29, No. 1, 504-519.

\bibitem {bryc2}Bryc, W. (2001) Stationary Markov chains with linear
regressions. Stochastic Processes and Applications 93, 339-348.

\bibitem {bms}Bryc, W., Matysiak, W., Szab\l owski, P.J. (2003) Probabilistic
aspects of Al-Salam--Chihara polynomials. Proceedings of the American Mathematical Society (to appear). \\
WWW:\ \texttt{http://arXiv.org/abs/math/0304155}.

\bibitem {doob}Doob, J. (1953) Stochastic Processes. Wiley, New York.

\bibitem {matszab}Matysiak, W., Szab\l owski, P.J. (2002) A few remarks on
Bryc's paper on random fields with linear regressions. Annals of Probability
30, No. 3, 1486-1491.
\end{thebibliography}
\end{document}